\documentclass[journal]{IEEEtran}

\usepackage{cite}      

\usepackage[dvips]{graphicx}  

\usepackage{psfrag}    

\usepackage{url}       

\usepackage{amsfonts}
\usepackage{amsmath,amssymb}   
\usepackage{array}

\usepackage[latin1]{inputenc}
\usepackage[T1]{fontenc}
\usepackage[babel,final]{microtype} 

\DeclareMathOperator*{\argmin}{arg\,min}
\DeclareMathOperator{\I}{I}

\begin{document}

\newcommand{\tr}[1]{\ensuremath{\textnormal{trace}\left[{#1}\right]}}
\newcommand{\trace}{\ensuremath{\textnormal{trace}}	}
\newcommand{\veco}[1]{\ensuremath{\textnormal{vec}\left[{#1}\right]}}
\newcommand{\by}[2]{\ensuremath{\left[{#1}\times{#2}\right]}}
\newcommand{\kron}[2]{\ensuremath{\left({#1}\otimes{#2}\right)}}
\newcommand{\kronnp}[2]{\ensuremath{{#1}\otimes{#2}}} 
\newcommand{\abs}[1]{\ensuremath{\left|#1\right|}}
\newcommand{\norm}[1]{\ensuremath{\left|\left|#1\right|\right|}}
\newcommand{\expectation}[1]{\ensuremath{\mathbb{E}\left[#1\right]}}

\title{Mathematically Equivalent Approaches for Equality Constrained Kalman Filtering}

\author{Nachi~Gupta
\thanks{The author is with the Oxford University Computing Laboratory, Numerical Analysis Group, Wolfson Building, Parks Road, Oxford OX1 3QD, U.K., and can be reached via e-mail at \texttt{nachi@comlab.ox.ac.uk}.}}
\maketitle

\begin{abstract}
Kalman Filtering problems often have inherent and known constraints in the physical dynamics that are not exploited despite potentially significant gains (e.g., fixed speed of a motor).  In this paper, we review existing methods and propose some new ideas for filtering in the presence of equality constraints.  We then show that three methods for incorporating state space equality constraints are mathematically equivalent to the more general ``Projection'' method, which allows different weighting matrices when projecting the estimate.   Still, the different approaches have advantages in implementations that may make one better suited than another for a given application.
\end{abstract}



%
\IEEEpeerreviewmaketitle

\section{Introduction}
%
%
%
%

The Kalman Filter is the optimal estimator for dynamical systems with white process noise and measurement noise.  Since the inception of the Kalman Filter in 1960, a vast amount of research has gone into different extensions -- for example, to allow for nonlinear systems~\cite{JUD1995,KB2000,E1994,HM1998,E2006}, non-Gaussian noise distributions~\cite{Haug2005,AMGCS2002,K1996,GSS1993}, better numerical stability~\cite{VW2001,RGYU1999,VV1986,C1993,WLM1992,TB1978,TB1980,BT1977}, and state space constraints~\cite{SC2002,SAP1988,TS1988,WD1992,CWCP2002,WCC2002,SS2003,SS2005,Q1989,GH2007, Gupta2007,WXZWL2005,AB1993,CSC2006,YM2003,YZFM1999,GBS1997,P1988,T1985,CC1993,JL2007,MMB1995,CCC1991,D1992,SSKM1998,KB2007,YB2006,TCTAB2007}.   Incorporating these extensions gives rise to many sub-fields of Kalman Filtering. In the case of handling nonlinearities in the underlying system, there are a number of proposed models capturing different amounts of detail. This paper will focus on the last problem of incorporating state space constraints (that is, improving the best estimate given by the filtration process by accounting for known impossibilities). This is a small sub-field of Kalman Filtering, which has become more popular in just the past few decades and is growing rapidly. Specifically, we will focus on equality constraints.

We discuss a few distinct approaches to generalizing an equality constrained Kalman Filter.  The first approach is to augment the measurement space of the Kalman Filter with the equality constraints as noise-free measurements (also called pseudo-measurements) \cite{TS1988, SAP1988, TS1990}.  The second approach is to find the unconstrained estimate from a Kalman Filter and project it down to the equality constrained space~\cite{WD1992,SC2002}.  The third approach is to restrict the optimal Kalman Gain so the updated state estimate lies in the constrained space.  The fourth approach is to fuse the state prediction with the measurement in the constrained space.  The second method ends up being a generalization of the other three methods, i.e., the other methods are all special cases of the second method.  Proofs for this will be given (see also~\cite{Gupta2007,G2008,GH2007}).  Yet another approach to this problem would be to reduce the state space by the dimension of the constraints (i.e., to introduce an explicit coordinate system in the constrained space).  This leads to a state space that might not have an intuitive meaning in terms of the propagation equations.  This approach, while valid, is not discussed here.

Analogous to the way a Kalman Filter can be extended to solve problems containing nonlinearities, linear equality constrained filtering can be extended to problems with nonlinear constraints by linearizing locally (or another scheme motivated by how nonlinear filters handle the nonlinearities).  The accuracy achieved by methods dealing with nonlinear constraints will naturally depend on the structure and curvature of the nonlinear function itself.  

\section{Kalman Filter} \label{sec::kf}

The Kalman Filter is a formulation of the recursive least squares algorithm, which makes only one pass through the data such that it can wait for each measurement to come in real time and make an estimate at that time given all the information from the past.  In addition, the Kalman Filter holds a minimal amount of information in memory at each time for a cheap computational cost in solving the optimization problem.  A discrete-time Kalman Filter attempts to find the best running estimate for a recursive system governed by the following model.\footnote{Dealing with noise that is not normally distributed doesn't lend itself well to the framework of the Kalman Filter; for an arbitrary distribution, a Particle Filter~\cite{RA2004} could be used, and for noises that have heavy tail distributions such as power laws and L\'evy laws, the ``Kalman-L\'evy'' Filter has been proposed~\cite{SI2000,GPR2003}.}

\begin{equation} \label{kfsm} 
x_{k} = F_{k,k-1} x_{k-1} + u_{k,k-1}, \qquad u_{k,k-1} \sim \mathcal{N}\left(0,Q_{k,k-1}\right) 
\end{equation}

\begin{equation} \label{kfmm} 
z_{k} = H_{k} x_{k} + v_{k}, \qquad v_{k} \sim \mathcal{N}\left(0,R_{k}\right) 
\end{equation}

Here $x_{k}$ is an $n$-vector that represents the true state of the underlying system\footnote{The subscript $k$ means for the $k$-th time step, and all vectors in this paper are column vectors (unless of course we are taking the transpose of the vector).} and $F_{k,k-1}$ is an $n \times n$ matrix that describes the transition dynamics of the system from  $x_{k-1}$ to $x_{k}$.  The measurement made by the observer is an $m$-vector $z_{k}$, and $H_{k}$ is an $m \times n$ matrix that transforms a vector from the state space into the appropriate vector in the measurement space.  The noise terms $u_{k,k-1}$ (an $n$-vector) and $v_{k}$ (an $m$-vector) encompass errors in $F_{k,k-1}$ and $H_{k}$ and are normally distributed with mean 0 and covariances given by $n \times n$ matrix $Q_{k,k-1}$ and $m \times m$ matrix $R_{k}$, respectively.  At each iteration, the Kalman Filter makes a state prediction for $x_k$, denoted $\hat{x}_{k|k-1}$.  We use the notation ${k|k-1}$ since we will only use measurements provided until time-step $k-1$ in order to make the prediction at time-step $k$.  The state prediction error $\tilde{x}_{k|k-1}$ is defined as the difference between the true state and the state prediction, as below.

\begin{equation} \label{se1}
\tilde{x}_{k|k-1} = x_{k} - \hat{x}_{k|k-1}
\end{equation}

The covariance structure for the expected error on the state prediction is defined as the expectation of the outer product of the state prediction error.  We call this covariance structure the error covariance prediction and denote it $P_{k|k-1}$.

\begin{equation} \label{P-outer1}
P_{k|k-1} = \mathbb{E}\left[\left(\tilde{x}_{k|k-1}\right)\left(\tilde{x}_{k|k-1}\right)'\right]
\end{equation}

The filter will also provide an updated state estimate for $x_{k}$, given all the measurements provided up to and including time step $k$.  We denote these estimates by $\hat{x}_{k|k}$.  We similarly define the state estimate error $\tilde{x}_{k|k}$ as below.

\begin{equation} \label{se2}
\tilde{x}_{k|k} = x_{k} - \hat{x}_{k|k}
\end{equation}

The expectation of the outer product of the state estimate error represents the covariance structure of the expected errors on the state estimate, which we call the updated error covariance and denote $P_{k|k}$.

\begin{equation} \label{P-outer2}
P_{k|k} = \mathbb{E}\left[\left(\tilde{x}_{k|k}\right)\left(\tilde{x}_{k|k}\right)'\right]
\end{equation}

At time-step $k$, we can make a prediction for the underlying state of the system by allowing the state to transition forward using our model for the dynamics and noting that $\mathbb{E}\left[u_{k,k-1}\right] = 0$.  This serves as our state prediction.

\begin{equation} \label{kfsp} 
\hat{x}_{k|k-1} = F_{k,k-1} \hat{x}_{k-1|k-1} 
\end{equation}

If we expand the expectation in Equation \eqref{P-outer1}, we have the following equation for the error covariance prediction.\footnote{We use the prime notation on a vector or a matrix to denote its transpose throughout this paper.}

\begin{equation} \label{kfcp} 
P_{k|k-1} = F_{k,k-1} P_{k-1|k-1} F_{k,k-1}' + Q_{k,k-1}
\end{equation}

We can transform our state prediction into the measurement space, which is a prediction for the measurement we now expect to observe.

\begin{equation} \label{kfmp} 
\hat{z}_{k|k-1} = H_{k} \hat{x}_{k|k-1}
\end{equation}

The difference between the observed measurement and our predicted measurement is the measurement residual, which we are hoping to minimize in this algorithm.

\begin{equation} \label{kfi} 
\nu_{k} = z_{k} - \hat{z}_{k|k-1} 
\end{equation}

We can also calculate the associated covariance for the measurement residual, which is the expectation of the outer product of the measurement residual with itself, $\mathbb{E}\left[\nu_k \nu_k'\right]$.  We call this the measurement residual covariance.

\begin{equation} \label{kfic} 
S_{k} = H_{k} P_{k|k-1} H_{k}' + R_{k} 
\end{equation}

We can now define our updated state estimate as our prediction plus some perturbation, which is given by a weighting factor times the measurement residual.  The weighting factor, called the Kalman Gain, will be discussed below.

\begin{equation} \label{kfsu} 
\hat{x}_{k|k} = \hat{x}_{k|k-1} + K_{k}  \nu_{k}
\end{equation}

Naturally, we can also calculate the updated error covariance by expanding the outer product in Equation \eqref{P-outer2}.\footnote{The $\I$ in Equation \eqref{kfcu} represents the $n \times n$ identity matrix.  Throughout this paper, we use $\I$ to denote the same matrix, except in Appendix~\ref{app::kv}, in which $\I$ is the appropriately sized identity matrix.}

\begin{equation} \label{kfcu-JF} 
P_{k|k} = \left(\I - K_{k} H_{k}\right) P_{k|k-1}   \left(\I - K_{k} H_{k}\right)' + K_k R_k K_k'
\end{equation}

Now we would like to find the Kalman Gain $K_k$, which minimizes the mean square state estimate error, $\mathbb{E}\left[\norm{\tilde{x}_{k|k}}^2\right]$.  This is the same as minimizing the trace of the updated error covariance matrix above.\footnote{Note that $v'v = \tr{vv'}$ for any vector $v$.}  Expanding Equation \eqref{kfcu-JF}, we have the following.  After some calculus\footnote{The trace is minimized when the following matrix derivative is equal to zero: \\ $\frac{\partial\ \tr{P_{k|k}}}{\partial K_k} = -2 \left( H_k P_{k|k-1} \right)' + 2 K_k S_k = 0$.  Solving this for $K_k$ yields Equation \eqref{kfcu-JF}.}, we find the optimal gain that achieves this, written below.


\begin{equation} \label{kfkg} 
K_{k} = P_{k|k-1} H_{k}' S_{k}^{-1} 
\end{equation}

Substituting Equation \eqref{kfkg} into Equation \eqref{kfcu-JF} gives the following simplified form for the updated error covariance.  

\begin{equation} \label{kfcu}
P_{k|k} = \left(\I - K_{k} H_{k}\right) P_{k|k-1}
\end{equation}

In computation, one should avoid using this form and use Equation \eqref{kfcu-JF}, also called the Joseph Form.  While the Joseph Form requires more computation, it better preserves symmetry and reduces numerical loss of positive definiteness for the covariance matrix. 

The covariance matrices in the Kalman Filter provide us with a measure for uncertainty in our predictions and updated state estimate.  This is a very important feature for the various applications of filtering since we then know how much to trust our predictions and estimates.  Also, since the method is recursive, we need to provide an initial covariance that is large enough to contain the initial state estimate to ensure comprehensible performance.  




\subsection{Fusion Interpretation} \label{sec::fi}

We can also think of the Kalman Filter as a fusion of the state prediction with the measurement at each iteration.  Since we know the error covariance matrices for the state prediction and the measurement, we can take this fusion under a weighting and also calculate a covariance matrix for the best estimate.  Let us begin by re-writing our system in the following manner.\footnote{The superscript $F$ notation is used to denote the ``fusion'' filter.}

\begin{equation} \label{feq1} 
z_{k}^F = H_{k}^F x_{k} + v_{k}^F, \qquad v_{k}^F \sim \mathcal{N}\left(0,R_{k}^F\right) 
\end{equation}

Here $z_{k}^F$ and $v_{k}^F$ are augmented vectors, and $H_k^F$ is an augmented matrix (see Equations \eqref{feqz}, \eqref{feqH}, and \eqref{feqv}).  The first block of $z_{k}^F$ represents the prediction for the current time step, and the second block is the measurement. 

\begin{equation} \label{feqz}
z_{k}^F =  \begin{bmatrix}
	\hat{x}_{k|k-1} \\ 
	z_{k} 
\end{bmatrix}
\end{equation}

The matrix $H_{k}^F$ takes our state into the measurement space, as before.

\begin{equation} \label{feqH}
H_{k}^F  = \begin{bmatrix}
	\I \\
	H_{k}
\end{bmatrix}
\end{equation}



Now we define $v_{k}^F$ as the noise term, in which $v_{k}^F$ is normally distributed with mean 0 and covariance given by matrix matrix $R_{k}^F$.  

\begin{equation} \label{feqv}
v_{k}^F = \begin{bmatrix}
	\tilde{x}_{k|k-1} \\
	v_{k}
\end{bmatrix}
\end{equation}

The block diagonal elements of $R_{k}^F$ represent the covariance of each block of $v_{k}^F$.  Notice that $R_{k}^F$ contains no block off-diagonal elements implying no cross-correlations.  However, using this formulation, cross-correlations could be modelled easily.

\begin{equation} \label{feqR}
R_{k}^F = \begin{bmatrix}
	P_{k|k-1}  & 0 \\
	0 & R_{k} \\
\end{bmatrix}
\end{equation} 

This method of expressing our problem can be thought of as a fusion of the state prediction and the new measurement at each iteration.  The optimal estimate, as defined by the weighted least-squares method, for the system in Equation \eqref{feq1} is the minimizer of the cost function below.

\begin{equation}
J(x_k) = \left(z_{k}^F - H_{k}^F x_{k} \right)' \left(R_k^F\right)^{-1} \left( z_{k}^F - H_{k}^F x_{k} \right)
\end{equation}

The minimizer, which is found by standard calculus, is the least squares solution given below.

\begin{equation} \label{xFk}
\hat{x}_{k|k}^F = \left(\left(H_k^F\right)' \left(R_k^F\right)^{-1} H_k^F\right)^{-1} \left(H_k^F\right)' \left(R_k^F\right)^{-1} z_k^F
\end{equation}

The covariance for this solution is the following.

\begin{equation} \label{PFk-def}
P_{k|k}^F = \left(\left(H_k^F\right)' \left(R_k^F\right)^{-1} H_k^F\right)^{-1} 
\end{equation}

Some manipulation shows that this result is the same as that of the Kalman Filter.\footnote{For complete details on this derivation and extensions to nonlinear filtering, see~\cite{G2008}.}

\section{Incorporating Equality Constraints}

Equality constraints in this paper are defined as below, in which $A$ is a $q \times n$ matrix, $b$ a $q$-vector, and $x_k$, the state, is a $n$-vector, with $q \leq n $.\footnote{\label{fn::D}$A$ and $b$ can be different for different $k$.  We don't subscript each $A$ and $b$ to avoid confusion.  We assume these constraints are well defined throughout this paper -- i.e., no constraints conflict with one another to cause a null solution and no constraints are repeated.  More specifically, we assume $A$ has full row rank.  Note that under these conditions if $A$ was a square matrix, the constraints would completely determine the state.}

\begin{equation} \label{eq-constraints}
A x_k = b
\end{equation}

We would like our updated state estimate to satisfy the constraint at each iteration, as below.

\begin{equation} \label{kfsu-con}
A \hat{x}_{k|k} = b
\end{equation}

Similarly, we may also like the state prediction to be constrained, which would allow a better forecast for the system.\footnote{We do not discuss this point further here.  For more on this, please see~\cite{G2008}.}

\begin{equation}
A \hat{x}_{k|k-1} = b
\end{equation}

\subsection{Augmenting the Measurement Space} \label{sec::ams}

The first method that we discuss for incorporating equality constraints into a Kalman Filter is to ``observe'' the constraints at each iteration as noise-free measurements (or pseudo-measurements).  To illustrate this, we augment the linear constraints in Equations \eqref{eq-constraints} to the system shown in Equations \eqref{kfsm} and \eqref{kfmm} as measurements with zero variance.  Thus, we can re-write the system.\footnote{The superscript $A$ notation is used to denote the ``augmented'' constrained filter and bears no relation to the $A$ in Equation \eqref{eq-constraints}.  Also, note that the dimension of the state space hasn't changed (e.g., $x_k^A$ is the same size as $x_k$).}

\begin{equation} x_{k}^A = F_{k,k-1} x^A_{k-1} + u_{k,k-1}, \qquad u_{k,k-1} \sim \mathcal{N}\left(0,Q_{k,k-1}\right) \end{equation}

\begin{equation} z_{k}^A = H^A_{k} x^A_{k} + v^A_{k}, \qquad v_{k} \sim \mathcal{N}\left(0,R^A_{k}\right) \end{equation}

The next three equations show the construction of the augmentation in the measurement space.

\begin{equation}
z_k^A = 
\begin{bmatrix} 
	z_{k} \\
	b
\end{bmatrix}
\end{equation}

\begin{equation} \label{HkD}
H_k^A= 
\begin{bmatrix}
	H_{k} \\
	A 
\end{bmatrix}
\end{equation}

\begin{equation} \label{RD}
R_k^A=
\begin{bmatrix}
	R_{k} & 0 \\
	0 & 0
\end{bmatrix}
\end{equation}

The augmented state now forces $A x_k^A$ to be equal to $b$ exactly (i.e., with no noise term) at every iteration.\footnote{$x_k^A$ is still constructed in the same fashion as $x_k$.}  Let us now expand the equations for the Kalman Filter prediction and update to gain a stronger understanding of how the filter has changed.  

The state prediction from Equation \eqref{kfsp} becomes the following.

\begin{equation} \label{kfspd}
\hat{x}^A_{k|k-1} = F_{k,k-1} \hat{x}^A_{k-1|k-1}
\end{equation}

The error covariance prediction from Equation \eqref{kfcp} becomes the following.

\begin{equation} \label{kfcpd}
P^A_{k|k-1} = F_{k,k-1} P^A_{k-1|k-1} F_{k,k-1}' + Q_{k,k-1}
\end{equation}

The measurement prediction from Equation \eqref{kfmp} can then be written in the following form.

\begin{subequations} \label{kfmpd}
\begin{align} 
\hat{z}^A_{k|k-1} & = H^A_{k} \hat{x}^A_{k|k-1} \\
& =
	\begin{bmatrix}
		H_{k} \hat{x}^A_{k|k-1}\\
		A \hat{x}^A_{k|k-1}
	\end{bmatrix}
\end{align}
\end{subequations}

Similarly, we can express the measurement residual from Equation \eqref{kfi} in the following manner.

\begin{subequations} \label{kfid}
\begin{align} 
\nu^A_{k} & = z_{k}^A - \hat{z}^A_{k|k-1} \\
& =
	\begin{bmatrix}
		z_k - H_{k} \hat{x}^A_{k|k-1}\\
		b -A \hat{x}^A_{k|k-1}
	\end{bmatrix}
\end{align}
\end{subequations}

We expand the measurement residual covariance from Equation \eqref{kfic} below.

\begin{subequations} \label{kficd}
\begin{align} 
S^A_{k} & = H^A_{k} P^A_{k|k-1} \left(H^{D}_{k}\right)' + R^A_{k} \\
& = 
	\begin{bmatrix}
		H_{k} \\
		A 
	\end{bmatrix}
	P^A_{k|k-1}
	\begin{bmatrix}
		H_k' & A' 
	\end{bmatrix}
	+
	\begin{bmatrix}
		R_k & 0 \\
		0 & 0
	\end{bmatrix} \\
& = 
	\begin{bmatrix} 
		H_k P^A_{k|k-1} H_k' + R_k & H_k P^A_{k|k-1} A' \\
		A P^A_{k|k-1} H_k' & A P^A_{k|k-1} A'
	\end{bmatrix}
\end{align}
\end{subequations}

The Kalman Gain can now be written as below.

\begin{equation}
K^A_{k} = P^A_{k|k-1} \left(H^A_{k}\right)' \left(S^A_{k}\right)^{-1}
\end{equation}

In order to further expand this term, we denote $\left(S^A\right)^{-1}$ in the following block matrix form. 

\begin{equation} \label{Sinv}
\begin{bmatrix}
	\left(S_k^A\right)^{-1}_a & \left(S_k^A\right)^{-1}_b \\
	\left(S_k^A\right)^{-1}_c & \left(S_k^A\right)^{-1}_d
\end{bmatrix}
\end{equation}

We then expand the Kalman Gain in terms of the block structure of Equation \eqref{Sinv}.

\begin{subequations} \label{kfkgd}
\begin{align}
K_k^A & = P_{k|k-1}^A
	\begin{bmatrix}
		H_k' & A'
	\end{bmatrix}
	\begin{bmatrix}
		\left(S_k^A\right)^{-1}_a & \left(S_k^A\right)^{-1}_b \\
		\left(S_k^A\right)^{-1}_c & \left(S_k^A\right)^{-1}_d
	\end{bmatrix} \\
& = 
	\begin{bmatrix}
		P_{k|k-1}^A H_k' & P_{k|k-1}^A A' 
	\end{bmatrix}
	\begin{bmatrix}
		\left(S_k^A\right)^{-1}_a & \left(S_k^A\right)^{-1}_b \\
		\left(S_k^A\right)^{-1}_c & \left(S_k^A\right)^{-1}_d
	\end{bmatrix} \\
& = \label{kfkgd-split}
	\begin{bmatrix}
	\left(K_k^A\right)_a & \left(K_k^A\right)_b
	\end{bmatrix}
\end{align}
\end{subequations}

Here, we have used the following two terms to shorten the expression above.

\begin{subequations}
\begin{align}
\label{KDa} \left(K_k^A\right)_a & = P_{k|k-1}^A H_k' \left(S_k^A\right)^{-1}_a + P_{k|k-1}^A A' \left(S_k^A\right)^{-1}_c \\
\label{KDb} \left(K_k^A\right)_b & = P_{k|k-1}^A H_k' \left(S_k^A\right)^{-1}_b + P_{k|k-1}^A A' \left(S_k^A\right)^{-1}_d
\end{align}
\end{subequations}

Furthermore, the updated state estimate from Equation \eqref{kfsu} takes the following form.

\begin{equation}
\hat{x}^A_{k|k} = \hat{x}^A_{k|k-1} +K^A_{k}  \nu^A_{k} \label{kfsud}
\end{equation}

And the updated error covariance from Equation \eqref{kfcu} changes in the following way.

\begin{equation}
P^A_{k|k} = (\I - K^A_{k} H^A_{k}) P^A_{k|k-1} \label{kfcud}
\end{equation}

Methods using augmentation in Kalman Filters have appeared for different applications in the past (e.g., Fixed-Point Smoothing~\cite{GMS1988}, Bias Detection~\cite{Friedland1969}).  In order to gain a stronger understanding of the effects of augmentation in Kalman Filters, it can be helpful to read and understand these methods, as well -- though they are not relevant to equality constrained Kalman Filtering.

\subsubsection{Improvement gained over an Unconstrained Filter}

For a given iteration, we are interested in the improvement gained by using this method over a method that does not incorporate equality constraints.  In order to do so, we would like to find the constrained estimated $\hat{x}^A_{k|k}$ in terms of the unconstrained estimate $\hat{x}_{k|k}$ (and similarly the constrained error covariance matrix $P^A_{k|k}$ in terms of the unconstrained error covariance matrix $P_{k|k}$).  Suppose we start with the same previous estimate and error covariance matrix for both filters.

\begin{equation} \label{x0}
\hat{x}^A_{k-1|k-1} = \hat{x}_{k-1|k-1}
\end{equation}

\begin{equation} \label{P0}
P^A_{k-1|k-1} = P_{k-1|k-1}
\end{equation}

Thus, we consider the benefit of using the new constrained filter over the unconstrained Kalman Filter gained in one iteration.  We can re-write all the constrained filter's equations in terms of the corresponding equations of the unconstrained Kalman Filter.

Starting with Equation \eqref{kfspd}, we find that the state prediction remains the same over one iteration.

\begin{subequations}
\begin{align}
\hat{x}^A_{k|k-1} &\stackrel{\eqref{x0}}{=} F_{k,k-1} \hat{x}_{k-1|k-1} \\
&\stackrel{\eqref{kfsp}}{=} 
	\hat{x}_{k|k-1}
\end{align}
\end{subequations}

Similarly, we find the error covariance prediction from Equation \eqref{kfcpd} remains the same over one iteration.

\begin{subequations}
\begin{align}
P^A_{k|k-1} &\stackrel{\eqref{P0}}{=} F_{k,k-1} P_{k-1|k-1} F_{k,k-1}' + Q_{k,k-1} \\
&\stackrel{\eqref{kfcp}}{=} 
	P_{k|k-1}
\end{align}
\end{subequations}


The measurement prediction from Equation \eqref{kfmpd} is then modified as below.

\begin{subequations}
\begin{align} 
\hat{z}^A_{k|k-1} & \stackrel{\eqref{x0}}{=}
	\begin{bmatrix}
		H_{k} \hat{x}_{k|k-1}\\
		A \hat{x}_{k|k-1}
	\end{bmatrix} \\
& \stackrel{\eqref{kfmp}}{=}
	\begin{bmatrix}
		\hat{z}_{k|k-1}\\
		A \hat{x}_{k|k-1}
	\end{bmatrix}
\end{align}
\end{subequations}

For the measurement residual from Equation \eqref{kfid}, we arrive at the following.

\begin{subequations} \label{kfid2}
\begin{align} 
\nu^A_{k} & \stackrel{\eqref{x0}}{=}
	\begin{bmatrix}
		z_k - H_{k} \hat{x}_{k|k-1}\\
		b -A \hat{x}_{k|k-1}
	\end{bmatrix}\\
& \stackrel{\eqref{kfi}}{=} 
	\begin{bmatrix}
		\nu_k\\
		b -A \hat{x}_{k|k-1}
	\end{bmatrix}
\end{align}
\end{subequations}

The measurement residual covariance from Equation \eqref{kficd} can then be expressed as below.

\begin{subequations} \label{S_reduced}
\begin{align} 
S^A_{k} & \stackrel{\eqref{P0}}{=} 
	\begin{bmatrix}
		H_k P_{k|k-1} H_k' + R_k & H_k P_{k|k-1} A' \\
		A P_{k|k-1} H_k' & A P_{k|k-1} A'
	\end{bmatrix} \\
& \stackrel{\eqref{kfic}}{=} 
	\begin{bmatrix}
		S_k & H_k P_{k|k-1} A' \\
		A P_{k|k-1} H_k' & A P_{k|k-1} A'
	\end{bmatrix}
\end{align}
\end{subequations}

We are interested in finding $\left(S_k^A\right)^{-1}$ in a block structure.  We notice that is a saddle point matrix of the form given in Appendix~\ref{app::spm}.  The inverse of a saddle point matrix is given in a block matrix form in the appendix.  We can apply this to Equation \eqref{S_reduced}.\footnote{We know that $A_S$ as defined in Appendix~\ref{app::spm} will be nonsingular since it represents the measurement residual covariance $S_k$.  If this matrix was singular, this would mean there exists no uncertainty in our measurement prediction {\em or} in our measurement, and thus there would be no ability to filter.  Similarly, we know that $J_S$ as defined in Appendix~\ref{app::spm} must also be nonsingular, which is equal to $A P_{k|k-1} A'$ (see Equation \eqref{S_reduced}).  This term projects the predicted error covariance down to the constrained space.  For well defined constraints (see Footnote~\ref{fn::D}), this will never be singular -- it will have the same rank as $A$.}

\begin{subequations} \label{SDinv_a_simp}
\begin{align}
\left(S_k^A\right)^{-1}_a  \stackrel{\eqref{kfic}}{=} 
	& \nonumber\left(S_k\right)^{-1} + \left(S_k \right)^{-1} H_k P_{k|k-1} A' \\
	& \left(A P_{k|k-1} A' - A P_{k|k-1} H_k' \left(S_k \right)^{-1} H_k  P_{k|k-1} A'\right)^{-1} \\
	& A P_{k|k-1} H_k' \left(S_k \right)^{-1} \\
\stackrel{\eqref{DPD}}{=}
	&\left(S_k\right)^{-1} + \left(S_k \right)^{-1} H_k P_{k|k-1} A'  \left(A P_{k|k} A'\right)^{-1} \\
	& A P_{k|k-1} H_k' \left(S_k \right)^{-1} \\
\stackrel{\eqref{KDDPDDK}}{=}
	&\left(S_k\right)^{-1} + K_k' A' \left(A P_{k|k} A'\right)^{-1} A K_k
\end{align}
\end{subequations}

In a similar manner using Equations \eqref{kfic}, \eqref{DPD}, and \eqref{KDDPDDK}, we arrive at the following remaining terms in the block structure.

\begin{align} \label{SDinv_b_simp}
\left(S_k^A\right)^{-1}_b  = & - K_k' A' \left(A P_{k|k} A'\right)^{-1} 
\end{align}

\begin{align} \label{SDinv_c_simp}
\left(S_k^A\right)^{-1}_c  = & - \left(A P_{k|k} A' \right)^{-1} A K_k
\end{align}

\begin{align} \label{SDinv_d_simp}
\left(S_k^A\right)^{-1}_d  = & \left(A P_{k|k} A' \right)^{-1}
\end{align}


Applying this to Equations \eqref{KDa}, we can find the first part of the Kalman Gain.


\begin{subequations} \label{KDa_simp}
\begin{align}
\left(K_k^A\right)_a & \stackrel{\eqref{P0}}{=}  
	&& P_{k|k-1} H_k' \left(S_k^A\right)^{-1}_a + P_{k|k-1} A' \left(S_k^A\right)^{-1}_c \\
&\stackrel{\eqref{SDinv_a_simp},\eqref{SDinv_c_simp}}{=} 
	&& \nonumber P_{k|k-1} H_k' \left(S_k\right)^{-1}  \\
	&&& + P_{k|k-1} H_k' K_k' A' \left(A P_{k|k} A' \right)^{-1} A K_k \\
	&&&  - P_{k|k-1} A' \left(A P_{k|k} A' \right)^{-1} A K_k \\
&\stackrel{\eqref{kfkg}}{=} 
	&& K_k - \left(P_{k|k-1} - P_{k|k-1}  H_k' K_k' \right) \\
	&&& A' \left(A P_{k|k} A' \right)^{-1} A K_k \\
&\stackrel{\eqref{P-PHK}}{=} 
	&& K_k - P_{k|k} A' \left(A P_{k|k} A' \right)^{-1} A K_k 	
\end{align}
\end{subequations}

Following similar steps using Equations \eqref{P0}, \eqref{SDinv_b_simp}, \eqref{SDinv_d_simp}, and \eqref{P-PHK}, we can arrive at the other part of the Kalman Gain.

\begin{equation} \label{KDb_simp}
\left(K_k^A\right)_b = P_{k|k} A' \left(A P_{k|k} A'\right)^{-1}
\end{equation}

We can then substitute our expressions for $K^A_k$ directly into Equation \eqref{kfsud} to find a simplified form of the updated state estimate.

\begin{subequations} \label{ckfx}
\begin{align}
\hat{x}^A_{k|k} & \stackrel{\eqref{x0}}{=}  && \hat{x}_{k|k-1} +K^A_{k}  \nu^A_{k} \\
&\stackrel{\eqref{kfkgd},\eqref{kfid2}}{=}&&
	\hat{x}_{k|k-1} + \left(K_k^A\right)_a \nu_k + \left(K_k^A\right)_b \left(b -A \hat{x}_{k|k-1}\right) \\
&\stackrel{\eqref{KDa_simp},\eqref{KDb_simp}}{=}
	&& \nonumber\hat{x}_{k|k-1} + K_k \nu_k - P_{k|k} A' \left(A P_{k|k} A' \right)^{-1} A K_k \nu_k \\
	&&& + P_{k|k} A' \left(A P_{k|k} A'\right)^{-1} \left(b -A \hat{x}_{k|k-1} \right) \\
&\stackrel{\eqref{kfsu}}{=}&& \nonumber
	\hat{x}_{k|k} - P_{k|k} A' \left(A P_{k|k} A' \right)^{-1} A \left( \hat{x}_{k|k} - \hat{x}_{k|k-1}\right) \\
	&&& + P_{k|k} A' \left(A P_{k|k} A'\right)^{-1} \left(b -A \hat{x}_{k|k-1} \right) \\
&=&&
	\hat{x}_{k|k} - P_{k|k} A' \left(A P_{k|k} A' \right)^{-1} \left(A \hat{x}_{k|k} - b \right)
\end{align}
\end{subequations}

Similarly, we can expand the updated error covariance in Equation \eqref{kfcud}.

\begin{subequations} \label{ckfp}
\begin{align} 
\label{ckfcua}P^A_{k|k} & \stackrel{\eqref{P0}}{=}  && \left(\I - K^A_{k} H^A_{k}\right) P_{k|k-1} \\
& \stackrel{\eqref{kfkgd},\eqref{HkD}}{=} && 
	\left(\I -  \left(K_k^A\right)_a H_k - \left(K_k^A\right)_b D_{k}\right) P_{k|k-1} \\
& \stackrel{\eqref{KDa_simp},\eqref{KDb_simp}}{=} &&
	\nonumber \left(\I - K_k H_k + P_{k|k} A' \left(A P_{k|k} A' \right)^{-1}  A K_k H_k \right. \\
	&&& \left. - P_{k|k} A' \left(A P_{k|k} A'\right)^{-1} A \right) P_{k|k-1} \\
& = &&
	\left(\I - K_k H_k \right) P_{k|k-1} \\
	&&& -  P_{k|k} A' \left(A P_{k|k} A' \right)^{-1} A \left(\I - K_k H_k \right)  P_{k|k-1} \\
& \stackrel{\eqref{kfcu}}{=} && 
	P_{k|k} -  P_{k|k} A' \left(A P_{k|k} A' \right)^{-1} A P_{k|k}
\end{align}
\end{subequations}

Equations \eqref{ckfx} and \eqref{ckfp} give us the improvement gained over an unconstrained Kalman Filter in a single iteration of the augmentation approach to constrained Kalman Filtering.  We see that the covariance matrix can only get smaller since we are subtracting a positive semi-definite matrix from $P_{k|k}$ above.\footnote{\label{footnote-smaller}If $A$ and $B$ are covariance matrices, we say $B$ is smaller than $A$ if $A-B$ is positive semidefinite.}

\subsection{Projecting the Unconstrained Estimate} \label{sec::pue}

The second approach to equality constrained Kalman Filtering is to run an unconstrained Kalman Filter and to project the estimate down to the constrained space at each iteration.  We can then feed the new constrained estimate into the unconstrained Kalman Filter and continue this process.  Such a method can be described by the following minimization problem for a given time-step $k$, in which $\hat{x}_{k|k}^P$ is the constrained estimate,  $\hat{x}_{k|k}$ is the unconstrained estimate from the Kalman Filter equations, and $W_k$ is any positive definite symmetric weighting matrix.\footnote{The superscript $P$ notation is used to denote the ``projected'' constrained filter.}

\begin{equation} \label{eq-proj-problem}
\hat{x}_{k|k}^P = \arg\min_{x} \left\{\left(x - \hat{x}_{k|k} \right)' W_k \left(x - \hat{x}_{k|k} \right) : A x = b\right\}
\end{equation}

The best constrained estimate is then given below.

\begin{equation} \label{bce-xP}
\hat{x}_{k|k}^P = \hat{x}_{k|k} - W_k^{-1} A' \left( A W_k^{-1} A' \right)^{-1} \left(A \hat{x}_{k|k} - b \right)
\end{equation}

If we choose $W_k = P_{k|k}^{-1}$, we obtain the same solution as Equation \eqref{ckfx}.  This is not obvious considering the differing approaches.  The updated error covariance under this assumption will be the same as Equation \eqref{ckfp} since $\hat{x}_{k|k}^P =  \hat{x}_{k|k}^A$.  This choice of $W_k$ is the most natural since it best describes the uncertainty in the state.  One can also show that this choice leads to the smallest updated error covariance matrix $P_{k|k}^P$ (see e.g.,~\cite{SC2002}).\footnote{That is, this choice of $W_k$ makes $P_{k|k}^P$ smaller than any other choice of $W_k$ (see Footnote \ref{footnote-smaller}).}


In the more general case, we can still find the updated error covariance as a function of the unconstrained Kalman Filter's updated error covariance matrix as before.  First, let us define the matrix $\Upsilon$ below.\footnote{Note that $\Upsilon A$ is a projection matrix, as is $\left(\I - \Upsilon A\right)$, by definition.  If $A$ is poorly conditioned, we can use a QR factorization to avoid squaring the condition number.}

\begin{equation}
\Upsilon = W_k^{-1} A' \left(A W_k^{-1} A' \right)^{-1}
\end{equation}

Equation \eqref{bce-xP} can then be re-written as follows.

\begin{equation}
\hat{x}_{k|k}^P = \hat{x}_{k|k} - \Upsilon\left(A \hat{x}_{k|k} - b \right)
\end{equation}

We can find a reduced form for $x_k - \hat{x}_{k|k}^P$ as below.\footnote{Remember $A x_k - b = 0$.}

\begin{subequations}
\begin{align}
x_k - \hat{x}_{k|k}^P &= x_k - \hat{x}_{k|k} +\Upsilon \left(A \hat{x}_{k|k} - b - \left(A x_k - b \right)\right) \\
&=
	x_k - \hat{x}_{k|k} +\Upsilon \left(A \hat{x}_{k|k} - A x_k\right) \\
&=
	-\left(\I - \Upsilon A \right) \left(\hat{x}_{k|k} - x_k\right)
\end{align}
\end{subequations}

Using the definition of the error covariance matrix, we arrive at the following expression.

\begin{subequations} \label{bce-PP}
\begin{align}
P_{k|k}^P &= \mathbb{E}\left[\left(x_k - \hat{x}_{k|k}^P\right)\left(x_k - \hat{x}_{k|k}^P\right)'\right] \\
&= 
	\mathbb{E}\left[\left(\I - \Upsilon A \right) \left(\hat{x}_{k|k} - x_k\right) \left(\hat{x}_{k|k} - x_k\right)' \left(\I - \Upsilon A \right)'\right] \\
&= 
	\left(\I - \Upsilon A \right) P_{k|k} \left(\I - \Upsilon A \right)' \\
&= 
	P_{k|k} - \Upsilon A P_{k|k} - P_{k|k} A' \Upsilon' +  \Upsilon A P_{k|k} A' \Upsilon' \\
&= 
	P_{k|k} - \Upsilon A P_{k|k}
\end{align}
\end{subequations}

In the projection framework, two different filters can be constructed -- one with a feedback loop, and one without.  That is, the Kalman Filter can be run in real-time, and as a post-processing step, the unconstrained estimate and updated error covariance matrix can be reformulated in the constrained space; or alternatively, the constrained estimate and its associated updated error covariance matrix can be fed back into the system in real-time.  A large benefit of incorporating constraints can be realized in both techniques, though the feedback system should generally outperform the system without feedback.
%


\subsection{Restricting the optimal Kalman Gain} \label{sec::rokg}

The third approach to equality constrained Kalman Filtering is to expand the updated state estimate term in Equation \eqref{kfsu-con} using Equation \eqref{kfsu}.  

\begin{equation} 
A \left(  \hat{x}_{k|k-1} + K_{k}  \nu_{k} \right) = b
\end{equation}

Then we can choose a Kalman Gain $K_k^R$, that restricts the updated state estimate to be in the constrained space.\footnote{The superscript $R$ notation is used to denote the ``restricted kalman gain'' constrained filter.}  In the unconstrained case, we chose the optimal Kalman Gain $K_k$, by solving the minimization problem below which yields Equation \eqref{kfkg}.

\begin{equation}
\begin{aligned}
K_k = \argmin_{K \in \mathbb{R}^{n \times m}} \trace & \left[ \left(\I - K H_{k}\right) P_{k|k-1}   \left(\I - K H_{k}\right)' + K R_k K'\right]
\end{aligned}
\end{equation}

Now we seek the optimal $K_k^R$ that satisfies the constrained optimization problem written below for a given time-step $k$.

\begin{equation} \label{min-con-chap3}
\begin{split}
 K_k^R = \argmin_{K \in \mathbb{R}^{n \times m}}\ & \trace \left[ \left(\I - K H_{k}\right) P_{k|k-1}   \left(\I - K H_{k}\right)'  + K R_k K'\right] \\
 \textnormal{s.t. } & A \left(  \hat{x}_{k|k-1} + K  \nu_{k} \right) = b 
\end{split}
\end{equation}

We will solve this problem using the method of Lagrange Multipliers.  First, we take the steps below, using the vec notation (column stacking matrices so they appear as long vectors, see Appendix~\ref{app::kv}) to convert all appearances of $K$ in Equation \eqref{min-con-chap3} into long vectors.  Let us begin by expanding the following term.



\begin{subequations}
\begin{gather}
\nonumber\trace\left[\left(\I - K H_{k}\right) P_{k|k-1}   \left(\I - K H_{k}\right)' + K R_k K' \right] \qquad \qquad \qquad \qquad \qquad \qquad \qquad\\
\begin{aligned}
&\stackrel{\hphantom{\eqref{kfic}}}{=} && \trace \left[ P_{k|k-1}  - K H_{k} P_{k|k-1} - P_{k|k-1} H_{k}' K'  \right. \\
&&& \left.+  K H_{k} P_{k|k-1}  H_{k}' K' + K R_k K' \right] \\
&\stackrel{\eqref{kfic}}{=} &&\trace \left[ P_{k|k-1}  - K H_{k} P_{k|k-1} - P_{k|k-1} H_{k}' K' +  K S_k K' \right] \\
&\stackrel{\hphantom{\eqref{kfic}}}{=}\label{trace-separated} &&\trace \left[ P_{k|k-1} \right] - \trace \left[ K H_{k} P_{k|k-1} \right] \\
&&& - \trace \left[ P_{k|k-1} H_{k}' K' \right] + \trace \left[ K S_k K' \right]
\end{aligned}
\end{gather}
\end{subequations}

We now expand the last three terms in Equation \eqref{trace-separated} one at a time.\footnote{We use the symmetry of $P_{k|k-1}$ in Equation \eqref{KHP} and the symmetry of $S_k$ in Equation \eqref{KSK}.}

\begin{equation} \label{KHP}
\begin{aligned}
\trace \left[ K H_{k} P_{k|k-1} \right] 
&\stackrel{\eqref{tr-ab}}{=}
	\veco{\left(H_k P_{k|k-1}\right)'}' \veco{K} \\
&\stackrel{\hphantom{\eqref{tr-ab}}}{=}
	\veco{P_{k|k-1} H_k'}' \veco{K}
\end{aligned}
\end{equation}


\begin{equation}
\trace \left[ P_{k|k-1} H_{k}' K' \right]
\stackrel{\eqref{tr-ab}}{=} 
	\veco{K}' \veco{P_{k|k-1} H_k'}
\end{equation}





\begin{equation} \label{KSK}
\begin{aligned}
\trace \left[ K S_k K' \right] 
&\stackrel{\eqref{tr-ab}}{=}
	\veco{K}' \veco{K S_k} \\
&\stackrel{\eqref{vec-ab}}{=} 
	\veco{K}' \kron{S}{\I} \veco{K}
\end{aligned}
\end{equation}

Remembering that $\trace \left[ P_{k|k-1} \right]$ is constant, our objective function can be written as below.

\begin{equation}
\begin{aligned}
\veco{K}' \left(\I \otimes S_k \right) \veco{K'} &- \veco{P_{k|k-1} H_k'}' \veco{K}\\
&- \veco{K}' \veco{P_{k|k-1} H_k'}
\end{aligned}
\end{equation}

Using Equation \eqref{vec-abc} on the equality constraints, our minimization problem is the following.



\begin{equation}
\begin{split}
K_k^R = \argmin_{K \in \mathbb{R}^{n \times m}}& \ \veco{K}' \kron{S_k}{\I} \veco{K} \\
& - \veco{P_{k|k-1} H_k'}' \veco{K} \\
& - \veco{K}'  \veco{P_{k|k-1} H_k'} \\
\textnormal{s.t. } &  \left( \nu_{k}' \otimes A \right) \veco{K}  = b - A \hat{x}_{k|k-1}
\end{split}
\end{equation}



Further, we simplify this problem so the minimization problem has only one quadratic term.  We complete the square as follows.  We want to find the unknown variable $\mu$ which will cancel the linear term.  Let the quadratic term appear as follows.  Note that the non-``$\veco{K}$" term is dropped as it is irrelevant for the minimization problem.


\begin{equation}
\left(\veco{K} + \mu \right)' \kron{S_k}{\I} \left( \veco{K} + \mu \right)
\end{equation}

The linear term in the expansion above is the following.

\begin{equation}
\veco{K}'  \kron{S_k}{\I} \mu + \mu' \kron{S_k}{\I} \veco{K}
\end{equation}

So we require that the two equations below hold.

\begin{equation}
\begin{aligned}
\kron{S_k}{\I} \mu &= -\veco{P_{k|k-1} H_k'} \\
\mu' \kron{S_k}{\I} &= -\veco{P_{k|k-1} H_k'}'
\end{aligned}
\end{equation}

This leads to the following value for $\mu$.

\begin{equation}
\begin{aligned}
\mu 
&\stackrel{\eqref{kron-inv}}{=}
	 - \kron{S_k^{-1}}{\I} \veco{P_{k|k-1} H_k'} \\
&\stackrel{\eqref{vec-abc}}{=}
	-\veco{P_{k|k-1} H_k' S_k^{-1}} \\
&\stackrel{\eqref{kfkg}}{=}
	-\veco{K_k}
\end{aligned}
\end{equation}





Using Equation \eqref{vec-sum}, our quadratic term in the minimization problem becomes the following.

\begin{equation}
\left(\veco{K - K_k} \right)' \kron{S_k}{\I} \left( \veco{K - K_k} \right)
\end{equation}

Let $\ell = \veco{K - K_k}$.  Then our minimization problem becomes the following.

\begin{equation}
\begin{aligned}
K_k^R = \argmin_{\ell \in \mathbb{R}^{mn}} & \ \ell' \kron{S_k}{\I} \ell \\
\textnormal{s.t. }&  \left( \nu_{k}' \otimes A \right) \left(\ell + \veco{K_{k}}\right)  = b - A \hat{x}_{k|k-1}
\end{aligned}
\end{equation}

We can then re-write the constraint taking the $\veco{K_k}$ term to the other side as below.

\begin{equation}
\begin{aligned}
\left( \nu_{k}' \otimes A \right) \ell & \stackrel{\hphantom{\eqref{vec-abc}}}= &&b - A \hat{x}_{k|k-1} - \left( \nu_{k}' \otimes A \right) \veco{K_{k}} \\
& \stackrel{\eqref{vec-abc}}{=} &&b - A \hat{x}_{k|k-1} -\veco{A K_{k} \nu_k} \\
& \stackrel{\hphantom{\eqref{vec-abc}}}= &&b - A \hat{x}_{k|k-1}  - A K_{k} \nu_k \\
& \stackrel{\eqref{kfsu}}= && b - A \hat{x}_{k|k}
\end{aligned}
\end{equation}

This results in the following simplified form.

\begin{equation} \label{first-SDPT3}
\begin{aligned}
K_k^R = \argmin_{\ell \in \mathbb{R}^{mn}}&\ \ell' \kron{S_k}{\I} \ell \\
\textnormal{s.t. }&  \left( \nu_{k}' \otimes A \right) l  = b - A \hat{x}_{k|k}
\end{aligned}
\end{equation}

We form the Lagrangian $\mathcal{L}$, for which we introduce $q$ Lagrange Multipliers in vector $ \lambda = \left( \lambda_1, \lambda_2, \ldots, \lambda_q \right)'$



\begin{equation}
\begin{aligned}
\mathcal{L} = & \ell' \kron{S_k}{\I} \ell -  \lambda' \left[\left( \nu_{k}' \otimes A \right) \ell  - b + A \hat{x}_{k|k}\right]
\end{aligned}
\end{equation}

We take the partial derivative with respect to $\ell$.\footnote{We used the symmetry of $\kron{S_k}{\I}$ here.}

\begin{equation} \label{partial1}
\frac{\partial \mathcal{L}}{\partial \ell} = 2 \ell' \kron{S_k}{\I} - \lambda' \left( \nu_{k}' \otimes A \right)  \\
\end{equation}

Similarly we can take the partial derivative with respect to the vector $\lambda$.

\begin{equation}
\frac{\partial \mathcal{L}}{\partial \lambda}  = \left( \nu_{k}' \otimes A \right) \ell  - b + A \hat{x}_{k|k}
\end{equation}

When both of these derivatives are set equal to the appropriate size zero vector, we have the solution to the system.  Taking the transpose of Equation \eqref{partial1}, we can write this system as $Mn = p$ with the following block definitions for $M,n$, and $p$.

\begin{equation} \label{M-matrix}
M = \begin{bmatrix}
	2  \kronnp{S_k}{\I} & \nu_{k} \otimes A' \\
	 \nu_{k}' \otimes A & 0_{\by{q}{q}}
\end{bmatrix}
\end{equation}

\begin{equation} \label{n-vector}
n = \begin{bmatrix}
	\ell \\
	\lambda
\end{bmatrix}
\end{equation}

\begin{equation} \label{p-vector}
p = \begin{bmatrix}
	0_{\by{mn}{1}} \\
	b - A \hat{x}_{k|k}
\end{bmatrix}
\end{equation}

We solve this system for vector $n$ in Appendix~\ref{app::Mnp}.  The solution for $\ell$ is copied below.

\begin{equation}
\left(\left[S_k^{-1} \nu_k \left(\nu_{k}' S_k^{-1} \nu_k \right)^{-1}\right] \otimes \left[A' \left(A A' \right)^{-1} \right]\right) \left(b - A \hat{x}_{k|k}\right)
\end{equation}

Bearing in mind that $b - A \hat{x}_{k|k} = \veco{b - A \hat{x}_{k|k}}$, we can use Equation \eqref{vec-abc} to re-write $l$ as below.\footnote{Here we used the symmetry of $S_k^{-1}$ and $\left(\nu_{k}' S_k^{-1} \nu_k \right)^{-1}$.}

\begin{equation}
\veco{A' \left(A A' \right)^{-1}\left(b - A \hat{x}_{k|k} \right)  \left(\nu_{k}' S_k^{-1} \nu_k \right)^{-1} \nu_k' S_k^{-1}}
\end{equation}

The resulting matrix inside the vec operation is then an $n$ by $m$ matrix.  Remembering the definition for $l$, we notice that $K - K_k$ also results in an $n$ by $m$ matrix.  Since both of the components inside the vec operation result in matrices of the same size, we can safely remove the vec operation from both sides.  This results in the following optimal constrained Kalman Gain $K_k^R$.

\begin{equation}
K_k - A' \left(A A' \right)^{-1}\left(A \hat{x}_{k|k} - b \right)  \left(\nu_{k}' S_k^{-1} \nu_k \right)^{-1} \nu_k' S_k^{-1}
\end{equation}

If we now substitute this Kalman Gain into Equation \eqref{kfsu} to find the constrained updated state estimate, we end up with the following.

\begin{equation}
\hat{x}_{k|k}^R = \hat{x}_{k|k} - A' \left(A A' \right)^{-1}\left(A \hat{x}_{k|k}  - b \right) 
\end{equation}

This is of course equivalent to the result of Equation \eqref{bce-xP} with the weighting matrix $W_k$ chosen as the identity matrix.  The error covariance for this estimate is given by Equation \eqref{bce-PP}.\footnote{We can use the unconstrained or constrained Kalman Gain to find this error covariance matrix.
Since the constrained Kalman Gain is suboptimal for the unconstrained problem, before projecting onto the constrained space, the constrained covariance will be different from the unconstrained covariance.  However, the difference lies exactly in the space orthogonal to which the covariance is projected onto by  Equation \eqref{bce-PP}.  The proof is omitted for brevity.}

\subsection{Fusion Approach} \label{sec::ams-fi}

The fourth approach to equality constrained Kalman Filtering is to augment the constraints onto the system using the fusion interpretation to the Kalman Filter.  In this case, we would like to fuse our state prediction with our measurement in the constrained space.  Our system is then defined as below.\footnote{The superscript $C$ notation is used to denote the ``augmented-fusion'' constrained filter.}

\begin{equation} \label{fceq1}  z_{k}^C = h_{k}^C(x_{k}) + v_{k}^C, \qquad v_{k}^C \sim \mathcal{N}\left(0,R_{k}^C\right) \end{equation}

Here $z_{k}^C$, $h_{k}^C$, and $v_{k}^C$ are all vectors, each having three distinct parts.  The first part represents the prediction for the current time-step, the second part is the measurement, and the third part is the equality constraint.  $z_{k}^C$ effectively still represents the measurement, with the prediction treated as a ``pseudo-measurement" with its associated covariance.

\begin{equation} \label{fceqz}
z_{k}^C =  \begin{bmatrix}
	\hat{x}_{k|k-1} \\ 
	z_{k} \\
	b
\end{bmatrix}
\end{equation}

The matrix $H_{k}^C$ takes our state into the measurement space, as before.

\begin{equation} \label{fceqH}
H_{k}^C  = \begin{bmatrix}
	\I \\
	H_{k} \\
	A
\end{bmatrix}
\end{equation}



Now we define $v_{k}^C$ as the noise term, in which $v_{k}^C$ is normally distributed with mean 0 and covariance given by matrix $R_{k}^C$.  

\begin{equation} \label{fceqv}
v_{k}^C = \begin{bmatrix}
	\tilde{x}_{k|k-1} \\
	v_{k} \\
	0
\end{bmatrix}
\end{equation}

The block diagonal elements of covariance matrix $R_{k}^C$ represent the covariance of each element of $v_{k}^C$.  We define the covariance of the state estimate error at time-step $k$ as $P_{k|k}$.  Notice that $R_{k}^C$ contains no block off-diagonal elements implying no cross-correlations.  However, in this formulation, cross-correlations can be modelled.

\begin{equation} \label{fceqR}
R_{k}^C = \begin{bmatrix}
	P_{k|k-1}  & 0 & 0\\
	0 & R_{k} & 0\\
	0 & 0 & 0
\end{bmatrix}
\end{equation} 

This method of expressing our problem can be thought of as a fusion of the state prediction and the new measurement at each iteration.  The solution and covariance for the problem given in Equation~\ref{fceq1} is printed below.

\begin{equation} \label{xcFk}
\hat{x}_{k|k}^C = \left(\left(H_k^C\right)' \left(R_k^C\right)^{-1} H_k^C\right)^{-1} \left(H_k^C\right)' \left(R_k^C\right)^{-1} z_k^C
\end{equation}

\begin{equation} \label{PcFk-def}
P_{k|k}^C = \left(\left(H_k^C\right)' \left(R_k^C\right)^{-1} H_k^C\right)^{-1} 
\end{equation}

However, the matrix $R_k^C$ is positive semi-definite now, and therefore singular, so the inverse is not well defined.  Let us look at the inverse of the following saddle point matrix.  The bottom left block will correspond exactly to the right-hand side of Equation \eqref{xcFk}, and the bottom right block will correspond to the negated right-hand side of Equation \eqref{PcFk-def}.  This can be verified by making the proper substitutions using Appendix~\ref{app::spm}.

\begin{equation}
\begin{bmatrix}
R_k^C & H_k^C \\
\left(H_k^C\right)' & 0
\end{bmatrix}^{-1}
\end{equation}

The statement for the inverse of the saddle point matrix made in Appendix~\ref{app::spm} also holds in the case where all inverses are replaced by the Moore-Penrose pseudo-inverse~\cite{CM1979}.  Taking the pseudo-inverse, we can correctly express Equations \eqref{xcFk} and \eqref{PcFk-def} below.\footnote{The pseudo-inverse may not be required if the matrix is invertible, the conditions for which are given in~\cite{NWL1992}.  Further, if the matrix is invertible, the pseudo-inverse will be the true inverse.}

\begin{equation}
\hat{x}_{k|k}^C = \begin{bmatrix}
		0 & \I
	\end{bmatrix}
	\begin{bmatrix}
		R_k^C & H_k^C \\
		\left(H_k^C\right)' & 0
	\end{bmatrix}^+ \\
	\begin{bmatrix}
		\I \\
		0
	\end{bmatrix}
	 z_k^C
\end{equation}

\begin{equation}
P_{k|k}^C = -\begin{bmatrix}
		0 & \I
	\end{bmatrix}
	\begin{bmatrix}
		R_k^C & H_k^C \\
		\left(H_k^C\right)' & 0
	\end{bmatrix}^+ \\
	\begin{bmatrix}
		0 \\
		\I
	\end{bmatrix}
\end{equation}

We have already shown that the Fusion interpretation of the filter is identical to the Kalman Filter in Section~\ref{sec::fi}, this method is mathematically equivalent to the method in Section~\ref{sec::ams}, in which we have also augmented pseudo-measurements.

\section{Nonlinear Equality Constraints} \label{sec::nec}

Since the equality constraints that we model are often nonlinear, it is important to make an extension to nonlinear equality constrained Kalman Filtering for the four methods discussed thus far.  We replace the linear equality constraint on the state space by the following nonlinear constraint $a_k\left(x_k\right) = b$, in which $a_k\left(\cdot\right)$ is a vector-valued function.  The method based on augmenting the constraints presented in Sections~\ref{sec::ams} and~\ref{sec::ams-fi} is trivially extended by using an Extended Kalman Filter.

Incorporating nonlinear equality constraints into the methods described in Section~\ref{sec::pue} and Section~\ref{sec::rokg} requires a more explicit change.  If we linearize our constraint, $a_k\left(x_k\right) = b$, about the current state prediction $\hat{x}_{k|k-1}$, we have the following.

\begin{equation}
a\left(\hat{x}_{k|k-1}^P\right) + A \left(x_k - \hat{x}_{k|k-1}^P \right) \approx b
\end{equation}

Here $A$ is defined as the Jacobian of $a$ evaluated at $\hat{x}_{k|k-1}^P$, similar to before.  This indicates, then, that the nonlinear constraint we would like to model can be approximated by the following linear constraint.

\begin{equation} \label{puenl}
A x_k \approx b + A \hat{x}_{k|k-1}^P - a\left(\hat{x}_{k|k-1}^P\right)
\end{equation}

Then our projected state is given as in Section~\ref{sec::pue}, with $A$ defined as above, and $b$ replaced by the right hand side of Equation \eqref{puenl}.  This linearizations is mathematically equivalent to the linearization step taken by the Extended Kalman Filter when augmenting the constraints.  

We can again take an iterative method, such as the Iterated Extended Kalman Filter, which takes multiple iterations and linearization per time-step.  For the fusion implementation, an iterative algorithm is given in~\cite{WCC2002,CWCP2002}.

\section{Discussion of Methods}

Thus far, we have discussed four different methods for incorporating equality constraints into a Kalman Filter, and we have shown that three of these are mathematically equivalent to the projection method under the assumption that the weighting matrix $W_k$ is chosen appropriately.  As such, the projection method is a more general formulation.  On the other hand, the augmentation methods provide a trivial extension to {\em soft} equality constrained Kalman Filtering by increasing the noise for the constraints, which are normally zero.  In implementations, there are some subtle differences.  For instance, the augmentation methods requires a minimal adjustment to codes for an existing Kalman Filter or an Extended Kalman Filter -- that is, we can pass in the augmented matrices and get the constrained estimate.  This is especially advantageous for codes that use variations of the standard linear Kalman Filter (e.g., an Unscented Kalman Filter).  

There is another more transparent difference between these methods.  In implementations, we are bound to receive numerical round-off error.  While these methods can be mathematically equivalent, we will not see the exact same result.  The round-off error that causes the most trouble occurs when the updated error covariance matrices $P_{k|k}^A$ or $P_{k|k}^P$ lose symmetry or positive definiteness.  A way around this is to use the Joseph Form of the updated error covariance, which we discuss in more detail below.  The error covariance matrix calculation using the fusion method $P_{k|k}^C$ should maintain positive definiteness and symmetry quite well in implementations as is.

\subsection{Numerical Preservation of the Updated Error Covariance} \label{sec::npuec}

We would like to find a form of Equation \eqref{ckfp} that preserves symmetry and positive definiteness better.   Let us start with the Joseph Form of the updated error covariance matrix given in Equation \eqref{kfcu-JF}.

\begin{equation} \label{cJF}
P_{k|k}^A = \left(\I - K^A_{k} H^A_{k}\right) P_{k|k-1} \left(\I - K^A H^A \right)' + K^A R^A \left(K^A\right)'
\end{equation}

First, let us define the projection $\Gamma_k$ below.

\begin{equation} \label{Gamma}
\Gamma_k = \I -  P_{k|k} A' \left(A P_{k|k} A' \right)^{-1} A
\end{equation}

Then we see that Equation \eqref{ckfp} can be written using $\Gamma_k$ as follows

\begin{equation} \label{GP}
P^A_{k|k} = \Gamma_k P_{k|k}
\end{equation}

From this, the following easily follows using Equations \eqref{kfcu} and \eqref{ckfcua}.

\begin{equation} \label{ImKHG}
\I - K^A_{k} H^A_{k} = \Gamma_k \left(\I - K_{k} H_{k}\right)
\end{equation}

Equation \eqref{ImKHG} will help us in reducing the term to the left of the ``+'' sign in Equation \eqref{cJF}.  Let us focus on the right-side, $K^A R^A \left(K^A\right)'$, for the moment.

\begin{subequations} \label{KaRKa}
\begin{align} 
K^A R^A \left(K^A\right)' & \stackrel{\eqref{kfkgd-split},\eqref{RD}}{=} 
	\begin{bmatrix}
		\left(K_k^A\right)_a & \left(K_k^A\right)_b
	\end{bmatrix}
	\begin{bmatrix}
		R_{k} & 0 \\
		0 & 0
	\end{bmatrix}
	\begin{bmatrix}
		\left(K_k^A\right)_a' \\
		\left(K_k^A\right)_b'
	\end{bmatrix} \\
&\stackrel{\hphantom{\eqref{kfkgd-split},\eqref{RD}}}{=} 
	\left(K_k^A\right)_a R_k \left(K_k^A\right)_a'	
\end{align}
\end{subequations}

In terms of $\Gamma_k$, we find the following to also be true.

\begin{equation} \label{KaGK}
\left(K_k^A\right)_a = \Gamma_k K_k
\end{equation}

We are now ready to use Equation \eqref{cJF} to find a simplified form for the constrained updated error covariance.

\begin{subequations} \label{GPG}
\begin{align}
P_{k|k}^A 
&\stackrel{\eqref{KaRKa}}{=} && 
	\left(\I - K^A_{k} H^A_{k}\right) P_{k|k-1} \left(\I - K^A H^A \right)'  \\
	&&& + \left(K_k^A\right)_a R_k \left(K_k^A\right)_a' \\
&\stackrel{\eqref{ImKHG}\eqref{KaGK}}{=} &&
	\Gamma_k \left(\I - K_{k} H_{k}\right) P_{k|k-1} \left(\I - K_k H_k \right)' \Gamma_k'  \\
	&&& + \Gamma_k K_k R_k K_k' \Gamma_k' \\
&\stackrel{\hphantom{\eqref{ImKHG}\eqref{KaGK}}}= &&
	\Gamma_k \left[ \left(\I - K_{k} H_{k}\right) P_{k|k-1} \left(\I - K_k H_k \right)' \right. \\
	&&& \left. + K_k R_k K_k'\right] \Gamma_k' \\
&\stackrel{\hphantom{\eqref{ImKHG}\eqref{KaGK}}}=&&
	\Gamma_k P_{k|k} \Gamma_k'
\end{align}
\end{subequations}

To summarize, we can use Equation \eqref{kfcu} or \eqref{kfcu-JF} to find $P_{k|k}$, and we can use Equation \eqref{GP} or \eqref{GPG} to find $P_{k|k}^A$.  In practice, we should use Equations \eqref{kfcu-JF} and \eqref{GPG}, when applicable, in order to maintain numerical stability.

\section{Conclusions}

We have presented four approaches for incorporating state space equality constraints into a Kalman Filter and shown that three of them are special cases of the ``Projection'' method, which is a generalization that allows different weighting matrices when projecting the estimate.  However, either of the two augmentation methods may prove easier in implementations since we can use existing Kalman Filter codes with minimal modifications.  With the augmentation methods, we can also make a natural extension to incorporate {\em soft} equality constraints, in which we allow the constraint to be slightly blurred by adding a proportionate amount of noise to the bottom right block entry of $R_k^A$ (see Equation \eqref{RD}).  For experiments, please refer to~\cite{G2008}.

\begin{appendix}

\subsection{Inverse of a Saddle Point Matrix} \label{app::spm}

$M_S$ is a saddle point matrix if it has the block form below.\footnote{The subscript $S$ notation is used to differentiate these matrices from any matrices defined earlier.}

\begin{equation} \label{spm}
M_S =
	\begin{bmatrix}
		A_S & B_S' \\
		B_S & -C_S
	\end{bmatrix}
\end{equation}

In the case that $A_S$ is nonsingular and the Schur complement $J_S = -\left(C_S + B_S A_S^{-1} B_S'\right)$ is also nonsingular in the above equation, it is known that the inverse of this saddle point matrix can be expressed in the analytic block representation below (see e.g.,~\cite{BGL2005}).

\begin{equation}
M_S^{-1} =
	\begin{bmatrix}
		A_S^{-1} + A_S^{-1} B_S'  J_S^{-1} B_S A_S^{-1} & -A_S^{-1} B_S' J_S^{-1} \\
		-J_S^{-1} B_S A_S^{-1} & J_S^{-1}
	\end{bmatrix}
\end{equation}





\subsection{Some Identities}

The following are identities that will prove useful in some of the earlier derivations of Section~\ref{sec::ams}.  The matrices in these identities are used as defined in Sections~\ref{sec::kf} and~\ref{sec::ams}.

\subsubsection*{First Identity}

\begin{subequations} \label{DPD}
\begin{align}
& A  P_{k|k-1} A' - A P_{k|k-1} H_k' \left(S_k \right)^{-1} H_k P_{k|k-1} A' \\
&\stackrel{\eqref{kfkg}}{=} 
	A P_{k|k-1} A' - A K_k H_k P_{k|k-1} A' \\
&= 
	A \left(\I - K_k H_k \right) P_{k|k-1} A' \\
&\stackrel{\eqref{kfcu}}{=} 
	A P_{k|k} A'
\end{align}
\end{subequations}


\subsubsection*{Second Identity}

In the first step below, we make use of the symmetry of $P_{k|k-1}$ and $\left(S_k\right)^{-1}$.

\begin{subequations} \label{KDDPDDK}
\begin{align}
&\left(S_k \right)^{-1}H_k P_{k|k-1} A' \left(A P_{k|k} A'\right)^{-1} A P_{k|k-1} H_k' \left(S_k \right)^{-1} \\
& = 
	 \nonumber \left(P_{k|k-1} H_k' \left(S_k \right)^{-1} \right)' A' \left(A P_{k|k} A'\right)^{-1} A P_{k|k-1} H_k' \left(S_k \right)^{-1}  \\
&\stackrel{\eqref{kfkg}}{=}
	K_k' A' \left(A P_{k|k} A'\right)^{-1} A K_k 
\end{align}
\end{subequations}


\subsubsection*{Third Identity}

\begin{subequations}  \label{P-PHK}
\begin{align}
& P_{k|k-1}  - P_{k|k-1} H_k' K_k'\\
& = \label{idsym1}
	P_{k|k-1} \left(\I - H_k'K_k' \right) \\
& = \label{idsym2}
	\left(\I - K_k H_k \right) P_{k|k-1} \\
& \stackrel{\eqref{kfcu}}{=}
	P_{k|k}
\end{align}
\end{subequations}


Again, we have made use of the symmetry of $P_{k|k-1}$ between Equations \eqref{idsym1} and \eqref{idsym2}.

\subsection{Kron and Vec} \label{app::kv}

In this appendix, we provide some definitions used earlier in the chapter.  Given matrix $A \in \mathbb{R}^{ m \times n}$ and $B \in \mathbb{R}^{p \times q}$, we can define the right Kronecker product as below.\footnote{The indices $m,n,p$, and $q$ and all matrix definitions are independent of any used earlier.  Also, the subscript notation $a_{1,n}$ denotes the element in the first row and $n$-th column of $A$, and so forth.}

\begin{equation}
\left( A \otimes B \right) = \begin{bmatrix}
a_{1,1} B & \cdots & a_{1,n} B \\
\vdots & \ddots & \vdots \\
a_{m,1} B & \cdots & a_{m,n} B
\end{bmatrix}
\end{equation}

Given appropriately sized matrices $A, B, C,$ and $D$ such that all operations below are well-defined, we have the following equalities.

\begin{equation} \label{kron-trans}
\left( A \otimes B \right)' = \left( A' \otimes B' \right)
\end{equation}

\begin{equation} \label{kron-inv}
\left( A \otimes B \right) ^{-1} = \left( A^{-1} \otimes B^{-1} \right)
\end{equation}

\begin{equation} \label{kron-dist}
\left( A \otimes B \right) \left( C \otimes D \right) = \left( AC \otimes BD \right)
\end{equation}

We can also define the vectorization of an $\by{m}{n}$ matrix $A$, which is a linear transformation on a matrix that stacks the columns iteratively to form a long vector of size $\by{mn}{1}$, as below.

\begin{equation}
\veco{A} = \begin{bmatrix}
a_{1,1} \\
\vdots \\
a_{m,1} \\
a_{1,2} \\
\vdots \\
a_{m,2} \\
\vdots \\
a_{1,n} \\
\vdots \\
a_{m,n}
\end{bmatrix}
\end{equation}

Using the vec operator, we can state the trivial definition below.

\begin{equation} \label{vec-sum}
\veco{A+B} = \veco{A} + \veco{B}
\end{equation}

Combining the vec operator with the Kronecker product, we have the following.

\begin{equation} \label{vec-ab}
\veco{AB} = \kron{B'}{\I} \veco{A}
\end{equation}

\begin{equation} \label{vec-abc}
\veco{ABC} = \left(C' \otimes A \right) \veco{B}
\end{equation}

We can express the trace of a product of matrices as below.

\begin{equation} \label{tr-ab}
\tr{AB} = \veco{B'}'\veco{A}
\end{equation}

\begin{subequations}
\begin{align}
\tr{ABC} &= 
	\label{trace-1} \veco{B}' \left(\I \otimes C\right) \veco{A} \\
&= 
	\label{trace-2} \veco{A}' \left(\I \otimes B \right) \veco{C} \\
&=
	\label{trace-3} \veco{A}' \left(C \otimes \I \right) \veco{B}
\end{align}
\end{subequations}

For more information, please see~\cite{LT1985}.

\subsection{Solution to the system $Mn=p$} \label{app::Mnp}

Here we solve the system $Mn=p$ from Equations \eqref{M-matrix}, \eqref{n-vector}, and \eqref{p-vector}, re-stated below, for vector $n$.

\begin{equation} \label{Mnp}
\begin{bmatrix}
	2  \kronnp{S_k}{\I} &  \nu_{k} \otimes A' \\
	 \nu_{k}' \otimes A  & 0_{\by{q}{q}}
\end{bmatrix} \begin{bmatrix}
	\ell \\
	\lambda
\end{bmatrix} = \begin{bmatrix}
	0_{\by{mn}{1}} \\
	b - A \hat{x}_{k|k}
\end{bmatrix}
\end{equation}

$M$ is a saddle point matrix with the following equations to fit the block structure of Equation \eqref{spm}.\footnote{We use Equation \eqref{kron-trans} with $B_S'$ to arrive at the same term for $B_s$ in Equation \eqref{Mnp}.}

\begin{align}
A_S & =  2  \kronnp{S_k}{\I} \\
B_S & = \nu_{k}' \otimes A  \\
C_S & = 0_{\by{q}{q}}
\end{align}

We can calculate the term $A_S^{-1} B_S'$.

\begin{subequations}
\begin{align}
A_S^{-1} B_S' & = \left[ 2\kron{S_k}{\I}\right]^{-1} \left( \nu_{k}' \otimes A \right)'  \\
&\stackrel{\eqref{kron-trans}\eqref{kron-inv}}{=}  \frac{1}{2} \kron{S_k^{-1}}{\I} \left( \nu_{k} \otimes A' \right) \\
&\stackrel{\eqref{kron-dist}}{=}  \frac{1}{2} \left( S_k^{-1} \nu_k \right) \otimes A'
\end{align}
\end{subequations}

And as a result we have the following for $J_S$.

\begin{subequations}
\begin{align}
J_S & = - \frac{1}{2} \left( \nu_{k}' \otimes A \right) \left[ \left( S_k^{-1} \nu_k \right) \otimes A' \right] \\
&\stackrel{\eqref{kron-dist}}{=} - \frac{1}{2} \left( \nu_{k}' S_k^{-1} \nu_k \right) \otimes \left(A A' \right) 
\end{align}
\end{subequations}

$J_S^{-1}$ is then, as below.

\begin{subequations}
\begin{align}
J_S^{-1} & = -2 \left[ \left( \nu_{k}' S_k^{-1} \nu_k \right) \otimes \left( A A' \right)\right]^{-1} \\
&\stackrel{\eqref{kron-inv}}{=} -2 \left(\nu_{k}' S_k^{-1} \nu_k  \right)^{-1} \otimes \left(A A' \right)^{-1}
\end{align}
\end{subequations}

For the upper right block of $M^{-1}$, we then have the following expression.

\begin{subequations}
\begin{align}
A_S^{-1} B_S' J_S^{-1} &= \left[\left( S_k^{-1} \nu_k \right) \otimes A' \right] \left[\left(\nu_{k}' S_k^{-1} \nu_k \right)^{-1} \otimes \left(A A' \right)^{-1}\right] \\
&\stackrel{\eqref{kron-dist}}{=}  \left[S_k^{-1} \nu_k \left(\nu_{k}' S_k^{-1} \nu_k \right)^{-1}\right] \otimes \left[A' \left(A A' \right)^{-1} \right]
\end{align}
\end{subequations}

Since the first block element of $p$ is a vector of zeros, we can solve for $n$ to arrive at the following solution for $\ell$.

\begin{equation}
\left(\left[S_k^{-1} \nu_k \left(\nu_{k}' S_k^{-1} \nu_k \right)^{-1}\right] \otimes \left[A' \left(A A' \right)^{-1} \right]\right) \left(b - A \hat{x}_{k|k}\right) \\
\end{equation}

The vector of Lagrange Multipliers $\lambda$ is given below.

\begin{equation}
-2 \left[\left(\nu_{k}' S_k^{-1} \nu_k \right)^{-1} \otimes \left(A A' \right)^{-1} \right] \left(b - A \hat{x}_{k|k}\right) 
\end{equation}

\end{appendix}

\bibliography{physreve}

\end{document}